\documentclass{daj}

\dajAUTHORdetails{%
  title = {A very sharp threshold for first order logic distinguishability of random graphs}, 
  author = {Itai Benjamini and Maksim Zhukovskii},
  plaintextauthor = {Itai Benjamini, Maksim Zhukovskii},
  runningtitle = {A very sharp threshold for first order logic distinguishability of random graphs}, 
  keywords = {random graphs, first order logic, descriptive complexity, graph isomorphism},
}   

\dajEDITORdetails{%
   year={2025},
   number={6},
   received={3 August 2022},   
   revised={23 January 2025},    
   published={14 July 2025},  
   doi={10.19086/da.138190},       
}   

\usepackage{amsmath,amsthm,amssymb,amsfonts,hyperref}

\newtheorem{theorem}{Theorem}[section]

\newtheorem{claim}{Claim}[section]

\theoremstyle{remark}
\newtheorem{remark}{Remark}[section]

\begin{document}

\begin{frontmatter}[classification=text]

\title{A very sharp threshold for first order logic distinguishability of random graphs} 

\author[itai]{Itai Benjamini}
\author[maks]{Maksim Zhukovskii}

\begin{abstract}
In this paper we find an integer $h=h(n)$ such that the minimum number of variables of a first order sentence that distinguishes between two independent uniformly distributed random graphs of size $n$ with the asymptotically largest possible probability $\frac{1}{4}-o(1)$ belongs to $\{h,h+1,h+2,h+3\}$. We also prove that the minimum (random) $k$ such that two independent random graphs are distinguishable by a first order sentence with $k$ variables belongs to $\{h,h+1,h+2\}$ with probability $1-o(1)$.
\end{abstract}
\end{frontmatter}



\section{Introduction}
\label{sc:intro}

We say that $G_1,G_2$ are {\it $k$-distinguishable}, if there exists a first order (FO) sentence $\varphi$ with $k$ variables such that $G_1\models\varphi$ while $G_2\not\models\varphi$. In this paper we address the following question. Let $G_n^1,G_n^2$ be chosen uniformly at random and independently from the set of all graphs on $[n]:=\{1,\ldots,n\}$. What is the minimum $k$ such that they are $k$-distinguishable?

The classical FO 0-1 law~\cite{Fagin,Glebskii} says that, for every {\it fixed} sentence, its truth value on the random graph on $n$ vertices approaches either 0 or 1 as $n\to\infty$. Therefore, fixed FO sentences do not distinguish between random graphs with high probability. However, if we allow the number of variables to grow with $n$, then this is not true. In this paper, we find the minimum number of variables of a FO sentence that distinguishes between two independent random graphs (up to a constant additive term bounded by 4).

Note that the number of variables is an important measure of complexity of a FO sentence. In particular, the problem of deciding whether a graph on $n$ vertices satisfies a FO sentence with $k$ variables can be solved in time $O(n^k)$ (see~\cite[Proposition 6.6]{Libkin}).\\

If we fix a sequence of sentences $\varphi=\varphi(n)$ in advance, then due to independence, ${\sf P}(G_n^1\models\varphi,G_n^2\not\models\varphi)\leq\frac{1}{4}$. However, if we choose $\varphi$ depending on the graphs (i.e. $\varphi$ is random), then it may happen that the graphs are distinguishable with probability close to 1. Below we give very sharp answers to the following questions.
\begin{enumerate}
\item What is the minimum (not random) number of variables of a FO sentence $\varphi$ such that ${\sf P}(G_n^1\models\varphi,G_n^2\not\models\varphi)$ is bounded away from 0 (or even equals $1/4-o(1)$)?
\item What is the minimum (random) $k$ such that $G_n^1,G_n^2$ are $k$-distinguishable? 
\end{enumerate}

To get a naive lower bound in both questions (that turns out to be very sharp), let us recall that the graphs satisfying {\it $k$-extension axioms} (that we define below) are not $(k+1)$-distinguishable. The extension axiom $\varphi_k$ is the FO sentence with $k+1$ variables saying that, for every $k$-set and its partition $A\sqcup B$, there exists a vertex $x$ outside the set adjacent to all vertices in $A$ and no vertex from $B$. Due to the following claim (see~\cite[Lemma 12.7]{Libkin}), in order to prove that whp $G_n^1,G_n^2$ are not $k$-distinguishable, it is sufficient to prove that with probability $1-o(1)$ ({\it with high probability} or briefly {\it whp}) the uniformly sampled random graph $G_n$ on $[n]$ satisfies $\varphi_{k-1}$.

\begin{claim}
If $G_1,G_2$ satisfy $\varphi_k$, then there is no sentence with $k+1$ variables distinguishing between $G_1,G_2$.
\end{claim}

By the union bound
\begin{equation}
 {\sf P}(G_n\not\models\varphi_k)\leq{n\choose k}2^k\left(1-\frac{1}{2^k}\right)^{n-k}=e^{f^{\mathrm{ext}}(k)},
\label{eq:bound_extension}
\end{equation}
where 
$$
 f^{\mathrm{ext}}(k)=k\ln  n-k\ln k+k+k\ln 2-n2^{-k}+O(\ln k).
$$
Note that $f^{\mathrm{ext}}(k)>k(1+o(1))$ if $k>\log_2 n$ and, for all $k$, 
$$
f^{\mathrm{ext}}(k+1)-f^{\mathrm{ext}}(k)=\ln n+n2^{-k-1}-O(\ln k).
$$
Therefore, there is a unique $\hat k\in[1,n]$ such that $f^{\mathrm{ext}}(\hat k)=0$:
$$
\hat k=\log_2 n-2\log_2\ln n+\log_2\ln 2+o(1).
$$ 
Clearly, $e^{f^{\mathrm{ext}}(\lfloor\hat k\rfloor)}\leq 1+o(1)$ and $e^{f^{\mathrm{ext}}(\lfloor\hat k\rfloor-1)}=o(1)$. From this, we immediately get that whp $G_n^1,G_n^2$ are not $\lfloor\hat k\rfloor$-distinguishable. It turns out that the lower bound is sharp up to a constant additive term.\\

\begin{theorem}
Let $\hat k\in[1,n]$ be the unique root of the equation $f^{\mathrm{ext}}(k)=0$.
\begin{enumerate}
\item The minimum $k$ such that there exists a FO sentence $\varphi$ with $k$ variables satisfying ${\sf P}(G_n^1\models\varphi,G_n^2\not\models\varphi)=\frac{1}{4}-o(1)$ belongs to $\{\lfloor\hat k\rfloor+1,\ldots,\lfloor\hat k\rfloor+4\}$.
\item The minimum $k$ such that $G_n^1$ and $G_n^2$ are $k$-distinguishable belongs to $\{\lfloor\hat k\rfloor+1,\lfloor\hat k\rfloor+2,\lfloor\hat k\rfloor+3\}$ whp.
\end{enumerate}
\label{th:main}
\end{theorem}

\begin{remark}
We actually prove a stronger version of Theorem~\ref{th:main} since quantifier depths of sentences that we use to prove upper bounds, coincide with their numbers of variables. Thus, the minimum quantifier depth of a deterministic FO sentence that distinguishes between $G_n^1$ and $G_n^2$ with probability $\frac{1}{4}-o(1)$ belongs to $\{\lfloor\hat k\rfloor+1,\ldots,\lfloor\hat k\rfloor+4\}$; 2) the minimum quantifier depth of a random FO sentence that distinguishes between $G_n^1$ and $G_n^2$ whp belongs to $\{\lfloor\hat k\rfloor+1,\lfloor\hat k\rfloor+2,\lfloor\hat k\rfloor+3\}$, as well.
\end{remark}

The lower bound in both parts of Theorem~\ref{th:main} is already proven. Note that a weaker but asymptotically sharp upper bound $\log_2 n(1+o(1))$ --- it differs from the lower bound by $\Theta(\ln\ln n)$ --- in both parts of Theorem~\ref{th:main} follows from the upper bound for the minimum quantifier depth of a FO sentence that defines a random graph whp~\cite{KPSS}. Moreover, from~\cite[Theorem 12]{KPSS}, the bound $k\leq\hat k+O(1)$ for both quantities from Theorem~\ref{th:main} follows for {\it infinitely many $n$}. The reason of why the random graph whp can be defined by a FO sentence with quantifier depth $O(\log n)$ is that whp every two vertices have non-isomorphic neighborhoods, and so the upper bound is obtained by recursion, where for every vertex we write a sentence describing its neighborhood. Our proof strategy is absolutely different. We prove the first part of Theorem~\ref{th:main} in Section~\ref{proof:1}. In Section~\ref{proof:2} we show that the same approach as in the proof in Section~\ref{proof:1} can be used to prove part 2 of Theorem~\ref{th:main}. 

\section{Proof of Theorem~\ref{th:main} part 1}
\label{proof:1}

The desired sentence that distinguishes between two random graphs says that there exists a set $S$ that induces a specific subgraph with some additional restrictions on adjacencies between this set and all  the other vertices. Let us first note that we are forced to make this property monotone (for the reason explained further). For that, we make the restrictions on subgraphs that are induced by $S$ to be downward closed. In other words, we construct a sequence of families of graphs $\mathcal{F}_0\subset\mathcal{F}_1\subset\ldots$ on the vertex set $[k]$ such that $\mathcal{F}_i$ is downward closed and also closed under isomorphism. It can be easily done by adding at each step a full single isomorphism class such that its number of edges is the same as on the previous step if possible, and has one more edge otherwise (the initial family is just a singleton consisting of an empty graph).  For the random variable $X(i)$ that counts the number of sets $S$ that induce a graph from $\mathcal{F}_i$ and satisfy the conditions on the edges between $S$ and $[n]\setminus S$, we will prove that for some specific choice of $i$, ${\sf E}X(i)=\ln 2+o(1)$. Note that typically (since whp the random graph on $[k]$ is asymmetric), $|\mathcal{F}_i|-|\mathcal{F}_{i-1}|=k!$ (and this is the maximum possible value of the difference), so we should care about the jumps of expectation when switching from $i-1$ to $i$ (not to skip the desired asymptotics). After we prove that such an $i$ exists, we also prove that the second factorial moment of $X(i)$ equals $(\ln 2)^2+o(1)$, that together with some $\varphi$-mixing condition implies that ${\sf P}(X=0)\to\frac{1}{2}$ (see Section~\ref{sc:moments}). Note that the monotonicity is used to prove that the mixing condition is satisfied. Clearly, this immediately implies the desired probability bound. We will define our constructions in three different ways depending on the value of $n$ in Section~\ref{sc:expectation}.\\

\subsection{Expectation}
\label{sc:expectation}

Let $k$ be large and let $\mathcal{F}_k$ be a family of graphs on $[k]$ closed under isomorphism. Consider a uniformly distributed random graph $G_n$ on $[n]$, and let $X^{\mathcal{F}}_k$ be the total number of $k$-sets $S\subset[n]$ such that $G_n[S]$ (subgraph of $G_n$ induced on $S$) is isomorphic to a graph from $\mathcal{F}_k$,  and there is no vertex $x$ outside $S$ adjacent to all vertices in $S$. Note that for convenience and in order to making the considered events monotone, here we weaken the `extendibility' condition in the extension axioms (we consider only the partition $S=A\sqcup B$, where $A=S$, and $B=\emptyset$). Luckily it does not affect significantly the upper bound. We get
$$
 {\sf E}X^{\mathcal{F}}_k={n\choose k}\frac{|\mathcal{F}_k|}{2^{{k\choose 2}}}\left(1-\frac{1}{2^k}\right)^{n-k}=e^{f_{\mathcal{F}}(k)},
$$
$$ 
 f_{\mathcal{F}}(k)= k\ln  n-k\ln k+k+\ln|\mathcal{F}_k|-\frac{k(k-1)}{2}\ln 2-n2^{-k}+O(\ln k).
$$
Let $k_{\mathcal{F}}\in\mathbb{R}$ be such that $f_{\mathcal{F}}(k_{\mathcal{F}})=0$. \\

Let $F$ be a graph on $[k]$ with a trivial automorphism group (it is easy to see that such a graph exists for large enough $k$ --- in particular, the random graph has no symmetries whp~\cite{Symmetry}) and $\mathcal{F}^*_k$ contains all isomorphic copies of $F$, i.e. $|\mathcal{F}^*_k|=k!$. Let $\mathcal{F}^0_k\subset\mathcal{F}^1_k\subset\ldots$ be a sequence of graph families on $[k]$ considered in the beginning of the proof, i.e.
\begin{itemize}
\item $\mathcal{F}^0_k$ is a singleton set containing only an empty graph,
\item for every $i$, $\mathcal{F}^i_k\setminus\mathcal{F}^{i-1}_{k}$ is a set of all isomorphic copies of a single graph $F_i$ on $[k]$,
\item for every $i$, $\mathcal{F}^i_k$ is downward closed (i.e., if $H_1\subset H_2$ are graphs on $[k]$, and $H_2\in\mathcal{F}^i_k$, then $H_1\in\mathcal{F}^i_k$ as well).
\end{itemize}
Note that, for every $i$, $|\mathcal{F}^i_k|-|\mathcal{F}^{i-1}_k|\leq k!$ implying that ${\sf E}X^{\mathcal{F}^i}_k-{\sf E}X^{\mathcal{F}^{i-1}}_k\leq{\sf E}X^{\mathcal{F}^*}_k$. If we can prove that, for some specific $k$, ${\sf E}X^{\mathcal{F}^0}_k\to 0$ while, for large enough $i$, ${\sf E}X^{\mathcal{F}^i}_k\to \infty$, then we may hope that there exists $i$ such that ${\sf E}X^{\mathcal{F}^i}_k\to \ln 2$. To show this, it is sufficient to prove that ${\sf E}X^{\mathcal{F}^*}_k\to 0$.

Now we let $\mathcal{F}_k=\mathcal{F}^*_k$ and $k=k_{\mathcal{F}^*}$. 
 Recall that ${\sf E}X^{\mathcal{F}}_{\lfloor k \rfloor}\leq 1$ 
 and ${\sf E}X^{\mathcal{F}}_{\lceil k\rceil}\geq 1$. Let $\Delta:=k-\lfloor k\rfloor$. Then
$$
f_{\mathcal{F}}(k)-
f_{\mathcal{F}}(\lfloor k\rfloor)=k^2\ln 2\left(2^{\Delta}-1+o(1)\right).
$$
In order to be certain that ${\sf E}X^{\mathcal{F}}_{\lfloor k \rfloor}\to 0$, we need $2^{\Delta}-1$ to be bounded away from 0. Moreover, in order to have ${\sf E}X^{\mathcal{F}^i}_{\lfloor k \rfloor}\to\infty$ for some $i$, we need $2^{\Delta}-1$ not to be very large. The sufficient restriction $0.05<\Delta\leq 0.55$ is considered in the next section. For other values of $\Delta$, we might slightly modify the sentence (and this is the reason why we get concentration in 4 points instead of 2).

\subsubsection{$0.05<\Delta\leq 0.55$}

If $\Delta>0.05$, then ${\sf E}X^{\mathcal{F}}_{\lfloor k \rfloor}=o(1)$. If, at the same time, $\Delta\leq 0.55$, then $2^{\Delta}-1<\frac{1}{2}$. Therefore, 
\begin{equation}
 \frac{2^{{\lfloor k\rfloor\choose 2}}}{k!}{\sf E}X^{\mathcal{F}}_{\lfloor k\rfloor}\gg e^{f_{\mathcal{F}}(k)}=1
\label{eq:higher_expectation_infinite}
\end{equation} 
implying that ${\sf E}X_{\lfloor k\rfloor}^{\bold{F}}\to\infty$, where $\bold{F}_k$ is the set of {\bf all} graphs on $[k]$. Let us recall that, for every $i$,
$$
{\sf E}X^{\mathcal{F}^{i-1}}_{\lfloor k \rfloor}+{\sf E}X^{\mathcal{F}^0}_{\lfloor k \rfloor}\leq {\sf E}X^{\mathcal{F}^i}_{\lfloor k \rfloor}\leq {\sf E}X^{\mathcal{F}^{i-1}}_{\lfloor k \rfloor}+{\sf E}X_{\lfloor k\rfloor}^{\mathcal{F}^*}.
$$
Note also that ${\sf E}X^{\mathcal{F}^0}_{\lfloor k \rfloor}=\frac{1}{k!}{\sf E}X^{\mathcal{F}^*}_{\lfloor k \rfloor}$ implying that ${\sf E}X^{\mathcal{F}^0}_{\lfloor k \rfloor}=o(1)$. Due to (\ref{eq:higher_expectation_infinite}), we get that there exists $\mu=\mu(k)$ such that 
$$
{\sf E}X^{\mathcal{F}^{\mu}}_{\lfloor k \rfloor}=\ln 2+o(1)
$$ 
as needed. 
Note that $\{X^{\mathcal{F}^{\mu}}_{\lfloor k \rfloor}=0\}$ is expressed by a FO sentence with 
$$
\lfloor k\rfloor +1=\lfloor k-0.05\rfloor +1\leq\lfloor\hat k\rfloor+2
$$ 
variables. The latter inequality follows from the fact that 
$f_{\mathcal{F}^0}(\hat k+1.05)\to\infty$ implying that $k<\hat k+1.05$. We compute the second moment of $X^{\mathcal{F}^{\mu}}_{\lfloor k \rfloor}$ in Section~\ref{sc:momentsX}.\\

\subsubsection{$0.55<\Delta\leq 0.95$}
\label{sc:modified_extension}

If $0.55<\Delta\leq 0.95$, then our goal is to modify the definition of $X_k$ in a way such that  the factor $\left(1-\frac{1}{2^k}\right)^{n-k}$ becomes $\left(1-\frac{3}{\cdot 2^k}\right)^{n-k}$. Having that, we will get that $k_{\mathcal{F}}$ becomes roughly $k_{\mathcal{F}}+\log_2 3$, and $\Delta$ becomes roughly $\Delta+\log_2 3-2\leq \log_2 3-1.05<0.55$. \\

Let $Y^{\mathcal{F}}_k$ be the total number of tuples $(v_1,v_2,S)$, where $v_1\neq v_2\in[n]$, $S\in{[n]\setminus\{v_1,v_2\}\choose k}$ such that $G_n[S]$ is equal to a graph from $\mathcal{F}_k$, and there are no vertices $x$ outside $S\sqcup\{v_1,v_2\}$ adjacent to all vertices from $S$ and at least 1 vertex among $v_1,v_2$. We have
$$
 {\sf E}Y^{\mathcal{F}}_k={n\choose 2}{n-2\choose k}\frac{|\mathcal{F}_k|}{2^{{k\choose 2}}}\left(1-\frac{3}{2^{k+2}}\right)^{n-k-2}=e^{g_{\mathcal{F}}(k)},
$$
$$ 
 g_{\mathcal{F}}(k)= (k+2)\ln  n-k\ln k+k+\ln|\mathcal{F}_k|-\frac{k(k-1)}{2}\ln 2-3n2^{-k-2}+O(\ln k).
$$ 
Let $k'_{\mathcal{F}}\in\mathbb{R}$ be such that $g_{\mathcal{F}}(k'_{\mathcal{F}})=0$. Note that 
$$
k'_{\mathcal{F}}=k_{\mathcal{F}}+\log_2\frac{3}{4}+o(1)=k_{\mathcal{F}}+\log_2 3-2+o(1).
$$
As above, letting $\mathcal{F}_k=\mathcal{F}^*_k$ and $k'=k'_{\mathcal{F}^*}$, we get that 
$$
\Delta':=k'-\lfloor k'\rfloor=\Delta+\log_2\frac{3}{2}+o(1)-\left\lfloor\Delta+\log_2\frac{3}{2}+o(1)\right\rfloor\in(0.1,0.55)
$$
since $0.55<\Delta\leq 0.95$. Then ${\sf E}Y^{\mathcal{F}}_{\lfloor k'\rfloor}=o(1)$ and
$$
g_{\mathcal{F}}(k')-
g_{\mathcal{F}}(\lfloor k'\rfloor)=(k')^2\ln 2\left(2^{\Delta'}-1+o(1)\right)<(k')^2\ln 2\left(0.47+o(1)\right)
$$
implying that
$$
 \frac{2^{{\lfloor k'\rfloor\choose 2}}}{k'!}{\sf E}Y^{\mathcal{F}}_{\lfloor k'\rfloor}\gg e^{g_{\mathcal{F}}(k')}=1.
$$ 
As above, we get that there exists $\mu'=\mu'(k)$ and $\mathcal{F}^{\mu'}_k$ containing all isomorphic copies of $\mu'$ non-isomorphic graphs on $[k]$ such that
$$
{\sf E}Y^{\mathcal{F}^{\mu'}}_{\lfloor k'\rfloor}=\ln 2+o(1).
$$ 
There exists a FO sentence describing $\{Y^{\mathcal{F}^{\mu'}}_{\lfloor k' \rfloor}=0\}$ with 
$$
\lfloor k'\rfloor +3\leq\lfloor\hat k+0.6\rfloor+3\leq\lfloor \hat k \rfloor+4
$$ 
variables. The inequalities follow from the fact that 
$g_{\mathcal{F}^*}(\hat k+0.6)\to\infty$. We compute the second moment of $Y^{\mathcal{F}^{\mu'}}_{\lfloor k' \rfloor}$ in Section~\ref{sc:momentsY}.\\

\subsubsection{$\Delta\leq 0.05$ or $\Delta>0.95$}

Finally, let either $\Delta\leq 0.05$ or $\Delta>0.95$ and  $Z^{\mathcal{F}}_k$ be the total number of tuples $(v_1,v_2,v_3,S)$, where $v_1,v_2,v_3\in[n]$ are distinct, $S\in{[n]\setminus\{v_1,v_2,v_3\}\choose k}$ such that $G_n[S]$ is equal to a graph from $\mathcal{F}_k$, and there are no vertices $x$ outside $S\sqcup\{v_1,v_2,v_3\}$ adjacent to all vertices from $S$ and either to $v_3$ or to both $v_1,v_2$. We have
$$
 {\sf E}Z^{\mathcal{F}}_k=3{n\choose 3}{n-3\choose k}\frac{|\mathcal{F}_k|}{2^{{k\choose 2}}}\left(1-\frac{5}{2^{k+3}}\right)^{n-k-3}=e^{h_{\mathcal{F}}(k)},
$$
$$ 
 h_{\mathcal{F}}(k)= (k+3)\ln  n-k\ln k+k+\ln|\mathcal{F}_k|-\frac{k(k-1)}{2}\ln 2-5n2^{-k-3}+O(\ln k).
$$ 
As above, letting $k''_{\mathcal{F}}\in\mathbb{R}$ be such that $f_{\mathcal{F}}(k''_{\mathcal{F}})=0$, $\mathcal{F}_k=\mathcal{F}^*_k$, $k''=k''_{\mathcal{F}^*}$, $\Delta'':=k''-\lfloor k''\rfloor$, we get that
$$
\Delta'':=k''-\lfloor k''\rfloor=\Delta+\log_2\frac{5}{4}+o(1)-\left\lfloor\Delta+\log_2\frac{5}{4}+o(1)\right\rfloor\in(0.25,0.4)
$$
implying 
$$
{\sf E}Z^{\mathcal{F}}_{\lfloor k''\rfloor}=o(1),\quad
 \frac{2^{{\lfloor k''\rfloor\choose 2}}}{k''!}{\sf E}Z^{\mathcal{F}}_{\lfloor k''\rfloor}\gg e^{h_{\mathcal{F}}(k'')}=1
$$ 
and the existence of $\mu''=\mu''(k)$ and $\mathcal{F}^{\mu''}_k$ containing all isomorphic copies of $\mu''$ non-isomorphic graphs on $[k]$ such that
$$
{\sf E}Z^{\mathcal{F}^{\mu''}}_{\lfloor k''\rfloor}=\ln 2+o(1).
$$ 
There exists a FO sentence describing $\{Z^{\mathcal{F}^{\mu''}}_{\lfloor k'' \rfloor}=0\}$ with 
$$
\lfloor k''\rfloor +4=\lfloor k''-0.25\rfloor +4\leq\lfloor\hat k-0.09\rfloor+4\leq\lfloor \hat k \rfloor+4
$$ 
variables. The inequalities follow from the fact that 
$h_{\mathcal{F}^*}(\hat k+0.16)\to\infty$. We compute the second moment of $Z^{\mathcal{F}^{\mu''}}_{\lfloor k'' \rfloor}$ in Section~\ref{sc:momentsZ}.\\

\subsection{The second moment and $\varphi$-mixing}
\label{sc:moments}

Let us recall that the $m$th factorial moment of a random variable $\xi$ is ${\sf E}\xi^{(m)}={\sf E}\xi(\xi-1)\ldots(\xi-m+1)$ and that, for $\xi\sim\mathrm{Pois}(\lambda)$, ${\sf E}\xi^{(m)}=\lambda^m$. It is well known that convergences of moments imply the convergence in distribution
(in particular, in the case of the convergence to a Poisson random variable), see, e.g., \cite[Theorem 30.2]{Billingsley}. Therefore, in order to prove the first part of Theorem~\ref{th:main}, it is sufficient to show that the $m$th factorial moments of all the random variables that we considered in the previous section converge to $(\ln 2)^m$. It turns out that it is easy to prove that the second factorial moment approaches $(\ln 2)^2$ while computing $m$th moments for $m\geq 3$ is not so evident. So, instead of computing all moments, we bound the probability of non-existence from below and above using~\cite[Lemmas 4.1, 4.2]{IRZZ}. Let us recall these bounds.

Let $A_1,\ldots,A_N$ be arbitrary events, and $X=\sum_{i=1}^N I[A_i]$ count the number of events that occur. Consider an arbitrary directed graph $\mathcal{D}$ on the vertex set $[N]$, and define 
$$
\psi_1=\sum_{i=1}^N\sum_{(j,i)\in\mathcal{D}:j<i}{\sf P}(A_i\cap A_j),\quad
\psi_2=\sum_{i=1}^N\sum_{(j,i)\in\mathcal{D}:j<i}{\sf P}(A_i){\sf P}(A_j),
$$
$$
\varphi=\max_i|{\sf P}(\cup_{(j,i)\notin\mathcal{D}:j<i}A_j|A_i)-{\sf P}(\cup_{(j,i)\notin\mathcal{D}:j<i}A_j)|.
$$
In~\cite{IRZZ}, it is proven that if $\psi_1,\psi_2,\varphi=o(1)$, then ${\sf P}(X=0)=\prod_{i=1}^N(1-{\sf P}(A_i))+o(1)$. 

Let $X:=X^{\mathcal{F}^{\mu}}_{\lfloor k \rfloor}$, $Y:=Y^{\mathcal{F}^{\mu'}}_{\lfloor k' \rfloor}$, $Z:=Z^{\mathcal{F}^{\mu''}}_{\lfloor k'' \rfloor}$. Note that $X=\sum_{i\in\left[{n\choose k}\right]} X_i$, where $X_i=I[A_i^X]$ is Bernoulli random variable equals to 1 if the event $A_i^X$ holds; 
\begin{multline*}
A_i^X=\{i\text{th }k\text{-set }S_i\text{ induces a copy of a graph from }\mathcal{F}^{\mu}\text{ in }G_n,\\
\text{and there are no vertices outside }S_i\text{ adjacent to all vertices of }S_i\}.
\end{multline*}
In the same way, we define $Y=\sum_{i\in\left[{n\choose 2}{n-2\choose k}\right]} Y_i$, $Z=\sum_{i\in\left[{n\choose 3}{n-3\choose k}\right]} Z_i$, where $Y_i=I[A_i^Y]$, $Z_i=I[A_i^Z]$; 
\begin{multline*}
A_i^Y=\{i\text{th tuple }(v_1^i,v_2^i,S_i)\text{ is such that }S_i\text{ induces a graph from }\mathcal{F}^{\mu'}\text{ in }G_n,\text{ and}\\
\text{there are no }u\notin S_i\sqcup\{v_1^i,v_2^i\}\text{ adjacent to all vertices of }S_i\text{ and at least 1 vertex among }v_1^i,v_2^i\}.
\end{multline*}
\begin{multline*}
A_i^Z=\{(v_1^i,v_2^i,v_3^i,S_i)\text{ is such that }S_i\text{ induces a graph from }\mathcal{F}^{\mu''}\text{ in }G_n,\text{ and}\\
\text{ no }u\notin S_i\sqcup\{v_1^i,v_2^i,v_3^i\}\text{ is adjacent to all vertices of }S_i\text{ and either to }v_3^i\text{ or to both }v_1^i,v_2^i\}.\end{multline*}
We define $\mathcal{D}^X,\mathcal{D}^Y,\mathcal{D}^Z$ in the most naive way: $(i,j)\in E(\mathcal{D}^{(\cdot)})$, if the $i$th and $j$th tuples are distinct and overlap. For these tuples of events we define $(\psi_1^X,\psi_2^X,\varphi^X)$, $(\psi_1^Y,\psi_2^Y,\varphi^Y)$, $(\psi_1^Z,\psi_2^Z,\varphi^Z)$ as suggested above. If we show that all these parameters approach $0$, then we immediately get ${\sf P}(X=0)=e^{-{\sf E}X(1+o(1))}=\frac{1}{2}+o(1)$, and the same holds for $Y$ and $Z$. Let us first show that $\varphi^X,\varphi^Y,\varphi^Z=o(1)$.

\subsubsection{$\varphi$-mixing}

Note that ${\sf P}\left(\cup_{(j,i)\notin\mathcal{D}^X:j<i}A_j^X|A_i^X\right)-{\sf P}\left(\cup_{(j,i)\notin\mathcal{D}^X:j<i}A_j^X\right)$ is bounded from above by the probability that there exists a set $S_j$ such that 
\begin{itemize}
\item $(j,i)\notin\mathcal{D}^X$, $j<i$,
\item $S_j$ induces a copy of a graph from $\mathcal{F}^{\mu}$,
\item  there are no vertices outside $S_i\cup S_j$ adjacent to all vertices of $S_j$,
\item and there is a vertex from $S_i$ that is adjacent to all vertices from $S_j$.
\end{itemize}
Therefore, 
\begin{equation}
{\sf P}\left(\cup_{(j,i)\notin\mathcal{D}^X:j<i}A_j^X|A_i^X\right)-{\sf P}\left(\cup_{(j,i)\notin\mathcal{D}^X:j<i}A_j^X\right)\leq{\sf E}X(1-(1-2^{-k})^k)(1+o(1))\to 0.
\label{eq:phiX_above}
\end{equation}
On the other hand,
$$
{\sf P}\left(\cup_{(j,i)\notin\mathcal{D}^X:j<i}A_j^X|A_i^X\right)=\frac{{\sf P}\left(\left[\cup_{(j,i)\notin\mathcal{D}^X:j<i}A_j^X\right]\cap A_i^X\right)}{{\sf P}(A_i^X)}.
$$
By construction, the events $\cup_{(j,i)\notin\mathcal{D}^X:j<i}A_j^X$ and $A_i^X$ are monotone decreasing. Therefore, by FKG inequality~\cite[Theorem 2.12]{Janson},
$$
 {\sf P}\left(\left[\cup_{(j,i)\notin\mathcal{D}^X:j<i}A_j^X\right]\cap A_i^X\right)\geq{\sf P}\left(\cup_{(j,i)\notin\mathcal{D}^X:j<i}A_j^X\right){\sf P}(A_i^X)
$$
implying that
\begin{equation}
{\sf P}\left(\cup_{(j,i)\notin\mathcal{D}^X:j<i}A_j^X|A_i^X\right)\geq{\sf P}\left(\cup_{(j,i)\notin\mathcal{D}^X:j<i}A_j^X\right).
\label{eq:phiX_below}
\end{equation}
Combining (\ref{eq:phiX_above}) with (\ref{eq:phiX_below}), we get that $\varphi^X=o(1)$ as needed.

The bounds for $\varphi^Y$ and $\varphi^Z$ are the same: the lower bound is immediate since $A_i^Y,A_i^Z$ are decreasing events. The upper bounds are similar to (\ref{eq:phiX_above}). Indeed, ${\sf P}\left(\cup_{(j,i)\notin\mathcal{D}^Y:j<i}A_j^Y|A_i^Y\right)-{\sf P}\left(\cup_{(j,i)\notin\mathcal{D}^X:j<i}A_j^Y\right)$ is bounded from above by the probability that there exists a tuple $(v_j^1,v_j^2,S_j)$ such that 
\begin{itemize}
\item $(j,i)\notin\mathcal{D}^Y$, $j<i$,
\item $S_j$ induces a copy of a graph from $\mathcal{F}^{\mu'}$,
\item there are no vertices outside $\{v_j^1,v_j^2,v_i^1,v_i^2\}\sqcup S_i\sqcup S_j$ adjacent to all vertices of $S_j$ and at least 1 vertex among $v_j^1,v_j^2$,
\item there is a vertex from $S_i\sqcup\{v_i^1,v_i^2\}$ that is adjacent to all vertices from $S_j$ and at least 1 vertex among $v_j^1,v_j^2$.
\end{itemize}
Then
$$
{\sf P}\left(\cup_{(j,i)\notin\mathcal{D}^Y:j<i}A_j^Y|A_i^Y\right)-{\sf P}\left(\cup_{(j,i)\notin\mathcal{D}^Y:j<i}A_j^Y\right)\leq
{\sf E}Y\left(1-\left(1-3\cdot 2^{-k-2}\right)^{k+2}\right)(1+o(1))\to 0.
$$
The proof of the upper bound
$$
{\sf P}\left(\cup_{(j,i)\notin\mathcal{D}^Z:j<i}A_j^Z|A_i^Z\right)-{\sf P}\left(\cup_{(j,i)\notin\mathcal{D}^Z:j<i}A_j^Z\right)\leq
{\sf E}Z\left(1-\left(1-5\cdot 2^{-k-3}\right)^{k+3}\right)(1+o(1))\to 0
$$
is identical.

\subsubsection{The second moment of $X$}
\label{sc:momentsX}

We here compute the second moment of $X:=X^{\mathcal{F}^{\mu}}_{\lfloor k \rfloor}$ or, more precisely, show that $\psi_1^X=o(1)$ and $\psi_2^X=o(1)$. We have
\begin{align*}
 \psi_1^X &\leq \sum_{s=1}^{k-1}{n\choose k}{k\choose s}{n-k\choose k-s}\frac{|\mathcal{F}^{\mu}_k|^2}{2^{2{k\choose 2}-{s\choose 2}}}\left(1-\frac{1}{2^s}+\frac{1}{2^s}\left(1-\frac{1}{2^{k-s}}\right)^2\right)^{n-2k+s}\\
 &\leq{n\choose k}^2\frac{|\mathcal{F}^{\mu}_k|^2}{2^{2{k\choose 2}}}\sum_{s=1}^{k-1}\frac{k^s}{s!}\frac{k^s}{n^s}2^{{s\choose 2}}\left(1-\frac{2}{2^k}+\frac{1}{2^{2k-s}}\right)^{n-2k+s}\\
&\leq({\sf E}X)^2(1+o(1))\sum_{s=1}^{k-1}\left(\frac{ek^2}{sn}2^{\frac{s-1}{2}}\right)^s\exp\left[\frac{n}{2^{2k-s}}\right].
\end{align*}
There is $s_0=k(1-o(1))$ such that $\frac{n}{2^{2k-s}}=o(1)$ for all $s\leq s_0$. If $s_0<s\leq k-1$, then $\frac{n}{2^{2k-s}}\leq\frac{n}{2^{k+1}}=\frac{k^2\ln 2}{4}(1+o(1))$. Therefore,
\begin{align*}
 \frac{\psi_1^X}{({\sf E}X)^2} & \leq(1+o(1))\biggl(\sum_{s=1}^{s_0}\left(\frac{ek^2}{sn}2^{\frac{s-1}{2}}\right)^s\\
 &\hspace{2cm}+\sum_{s=s_0+1}^{k-1}\exp\left[\left(-s\ln n+\frac{s^2}{2}\ln 2\right)(1+o(1))+\frac{k^2\ln 2}{4}(1+o(1))\right]\biggr)\\
& \leq \sum_{s=1}^{s_0}\left(\frac{k^2}{2^{k/2}}\right)^s+\sum_{s=s_0+1}^{k-1}\exp\left[-\left(1+o(1)\right)\frac{k^2\ln 2}{2}+\frac{k^2\ln 2}{4}(1+o(1))\right]=o(1).
\end{align*}
Moreover,
\begin{align*}
 \psi_2^X &=\sum_{s=1}^{k-1}{n\choose k}{k\choose s}{n-k\choose k-s}\left(\frac{|\mathcal{F}^{\mu}_k|}{2^{{k\choose 2}}}\left(1-\frac{1}{2^k}\right)^{n-k}\right)^2\\
&=({\sf E}X)^2\sum_{s=1}^{k-1}\frac{{k\choose s}{n-k\choose k-s}}{{n\choose k}}
\leq({\sf E}X)^2\sum_{s=1}^{k-1}\left(\frac{ek^2}{sn}\right)^s=o(1).
\end{align*}

Note that, since $A_i^X$ are monotone decreasing, by FKG inequality, we have that ${\sf P}(A_i^X\cap A_j^X)\geq{\sf P}(A_i^X){\sf P}(A_j^X)$ for all $i,j$. Therefore, $\phi_2^X\leq\phi_1^X$, and the same holds for $Y$ and $Z$. Therefore, we are actually not obliged to compute $\phi_2$ --- it is sufficient to prove that $\phi_1=o(1)$. In the next two sections, we show that $\phi_1^Y=o(1)$ , $\phi_1^Z=o(1)$, and this finishes the proof of the first part of Theorem~\ref{th:main}.

\subsubsection{The second moment of $Y$}
\label{sc:momentsY}

Let $\mathcal{F}$ be the family of graphs on $[k+2]$ that is obtained from $\mathcal{F}^{\mu'}_k$ by transforming each $F\in\mathcal{F}^{\mu'}_k$ into all possible $2^{1+2k}$ graphs on $[k+2]$ containing $F$ as a subgraph.

We separately count contributions to $\psi_1^Y$ of pairs of overlapping tuples that have at most $k$ vertices in common, and pairs of tuples that have $k+1$ vertices in common. For any two tuples $(v_i^1,v_i^2,S_i)$ and $(v_j^1,v_j^2,S_j)$ that share $s\leq k$ vertices and a vertex $z\notin\{v_i^1,v_i^2,v_j^1,v_j^2\}\cup S_i\cup S_j$, the probability that $z$ is adjacent to all vertices of $S_i\cup S_j$, at least 1 vertex among $v_i^1,v_i^2$ and at least 1 vertex among $v_j^1,v_j^2$, is at most $\frac{9}{2^{2k+4-s}}$. Moreover, if $s=k+1$, then this probability can be bounded by $\frac{5}{2^{k+3}}$ since its maximum is achieved when $S_i=S_j$, $v_i^2=v_j^2$, and $v_i^1\neq v_j^1$. In this case, there are exactly 5 ways to draw edges from $z$ to $(v_i^1,v_i^2,S_i)$ and $\{v_i^1,v_i^2,v_j^1,v_j^2\}\cup S_i\cup S_j$ that satisfies the mentioned restrictions.

Therefore, we get
\begin{align*}
 \psi_1^Y &\leq \sum_{s=1}^{k}{n\choose k+2}{k+2\choose s}{n-k-2\choose k+2-s}{k+2\choose 2}^2\frac{|\mathcal{F}|^2}{2^{2{k+2\choose 2}-{s\choose 2}}}\left(1-2\cdot\frac{3}{2^{k+2}}+\frac{9}{2^{2k+4-s}}\right)^{n-2k-4+s}\\
 &+{n\choose k+2}(k+2)(n-k-2){k+2\choose 2}^2\frac{|\mathcal{F}|^2}{2^{2{k+2\choose 2}-{k+1\choose 2}}}\left(1-2\cdot\frac{3}{2^{k+2}}+\frac{5}{2^{k+3}}\right)^{n-k-3}\\
&\leq(({\sf E}Y)^2+o(1))\biggl(\sum_{s=1}^{k}\frac{{k+2\choose s}{n-k-2\choose k+2-s}}{{n\choose k+2}}2^{{s\choose 2}}\exp\left[\frac{9n}{2^{2k+4-s}}\right]\\
&\hspace{7cm}+\frac{(k+2)(n-k-2)}{{n\choose k+2}}2^{{k+1\choose 2}}\exp\left[\frac{5n}{2^{k+3}}\right]\biggr)\\ 
&\leq(({\sf E}Y)^2+o(1))\left(\sum_{s=1}^{k}\left(\frac{e(k+2)^2}{sn}2^{\frac{s-1}{2}}\right)^s\exp\left[\frac{9n}{2^{2k+4-s}}\right]+\frac{(k+2)!}{n^k}2^{{k+1\choose 2}}\exp\left[\frac{5n}{2^{k+3}}\right]\right).
\end{align*}
There is $s_0=k(1-o(1))$ such that 
$$
\frac{n}{2^{2k+4-s}}=o(1)\text{ for all }s\leq s_0.
$$
If $s_0<s\leq k$, then 
$$
\frac{9n}{2^{2k+4-s}}\leq\frac{3}{4}\cdot\frac{3n}{2^{k+2}}=\frac{3k^2\ln 2}{8}(1+o(1)).
$$
In the same way, 
$$
\frac{5n}{2^{k+3}}=\frac{5k^2\ln 2}{12}(1+o(1)).
$$
Therefore,
\begin{align*}
 \frac{\psi_1^Y}{({\sf E}Y)^2} & \leq(1+o(1))\sum_{s=1}^{s_0}\left(\frac{e(k+2)^2}{sn}2^{\frac{s-1}{2}}\right)^s\\
&\hspace{2cm}+ \exp\left[\left(-k\ln n+\frac{k^2}{2}\ln 2\right)(1+o(1))+\frac{5k^2\ln 2}{12}(1+o(1))\right]\\
 &\hspace{3cm}+\sum_{s=s_0+1}^{k}\exp\left[\left(-s\ln n+\frac{s^2}{2}\ln 2\right)(1+o(1))+\frac{3k^2\ln 2}{8}(1+o(1))\right]\\
& \leq (1+o(1))\sum_{s=1}^{s_0}\left(\frac{k^2}{2^{k/2}}\right)^s+\sum_{s=s_0+1}^{k}\exp\left[-\left(1+o(1)\right)\frac{k^2\ln 2}{8}\right]\\
&\hspace{8cm}+\exp\left[-\left(1+o(1)\right)\frac{k^2\ln 2}{12}\right]=o(1).
\end{align*}

\subsubsection{The second moment of $Z$}
\label{sc:momentsZ}

Let $\mathcal{F}$ be the family of graphs on $[k+3]$ that is obtained from $\mathcal{F}^{\mu''}_k$ by transforming each $F\in\mathcal{F}^{\mu''}_k$ into all possible $2^{3+3k}$ graphs on $[k+3]$ containing $F$ as a subgraph.

For any two tuples $(v_i^1,v_i^2,v_i^3,S_i)$ and $(v_j^1,v_j^2,v_j^3,S_j)$ that share $s\leq k$ vertices and a vertex $z\notin\{v_i^1,v_i^2,v_j^1,v_j^2\}\cup S_i\cup S_j$, the probability that 
$z$ is adjacent to all vertices of $S_i\cup S_j$, either to $v_i^3$, or to both $v_i^1,v_i^2$, and either to $v_j^3$, or to both $v_j^1,v_j^2$, is at most $\frac{25}{2^{2k+6-s}}$. Moreover, if $s=k+1$, then this probability can be bounded by $\frac{17}{2^{k+5}}$ since its maximum is achieved when $S_i=S_j$, $v_i^3=v_j^3$, and $v_i^1,v_i^2,v_j^1,v_j^2$ are all distinct. Finally, if $s=k+2$, then this probability can be bounded by $\frac{9}{2^{k+4}}$ since its maximum is achieved when $S_i=S_j$, $v_i^3=v_j^3$, $v_i^2=v_j^2$, and $v_i^1\neq v_j^1$.

We have
\begin{align*}
 \psi_1^Z &\leq \sum_{s=1}^{k}9{n\choose k+3}{k+3\choose s}{n-k-3\choose k+3-s}{k+3\choose 3}^2\frac{|\mathcal{F}|^2}{2^{2{k+3\choose 2}-{s\choose 2}}}\left(1-2\cdot\frac{5}{2^{k+3}}+\frac{25}{2^{2k+6-s}}\right)^{n-2k-6+s}\\
 &\hspace{0.7cm}+9{n\choose k+3}{k+3\choose k+1}{n-k-3\choose 2}{k+3\choose 3}^2\frac{|\mathcal{F}|^2}{2^{2{k+3\choose 2}-{k+1\choose 2}}}\left(1-2\cdot\frac{5}{2^{k+3}}+\frac{17}{2^{k+5}}\right)^{n-k-5}\\
 &\hspace{0.7cm}+9{n\choose k+3}(k+3)(n-k-3){k+3\choose 3}^2\frac{|\mathcal{F}|^2}{2^{2{k+3\choose 2}-{k+2\choose 2}}}\left(1-2\cdot\frac{5}{2^{k+3}}+\frac{9}{2^{k+4}}\right)^{n-k-4}
 \end{align*}
 \begin{align*}
  \quad &\leq(({\sf E}Z)^2+o(1))\biggl(\sum_{s=1}^{k}\frac{{k+3\choose s}{n-k-3\choose k+3-s}}{{n\choose k+3}}2^{{s\choose 2}}\exp\left[\frac{25n}{2^{2k+6-s}}\right]\\
&\hspace{2.7cm}+\frac{{k+3\choose 2}{n-k-3\choose 2}}{{n\choose k+3}}2^{{k+1\choose 2}}\exp\left[\frac{17n}{2^{k+5}}\right]+\frac{(k+3)(n-k-3)}{{n\choose k+3}}2^{{k+2\choose 2}}\exp\left[\frac{9n}{2^{k+4}}\right]\biggr)\\
  \end{align*}
 \begin{align*}
\quad&\leq(({\sf E}Z)^2+o(1))\biggl(\sum_{s=1}^{k}\left(\frac{e(k+3)^2}{sn}2^{\frac{s-1}{2}}\right)^s\exp\left[\frac{25n}{2^{2k+6-s}}\right]\\
&\hspace{5cm}+\frac{(k+3)!}{n^k}2^{{k+1\choose 2}}\exp\left[\frac{17n}{2^{k+5}}\right]+\frac{(k+3)!}{n^{k+1}}2^{{k+2\choose 2}}\exp\left[\frac{9n}{2^{k+4}}\right]\biggr).
\end{align*}
As above we find $s_0=k(1-o(1))$ such that $\frac{n}{2^{2k+6-s}}=o(1)$ for all $s\leq s_0$. Then
\begin{align*}
 \frac{\psi_1^Z}{({\sf E}Z)^2} & \leq(1+o(1))\sum_{s=1}^{s_0}\left(\frac{e(k+3)^2}{sn}2^{\frac{s-1}{2}}\right)^s+\exp\left[\left(-k\ln n+\frac{k^2}{2}\ln 2+\frac{9k^2\ln 2}{20}\right)(1+o(1))\right]\\
&\hspace{2.5cm}+ \exp\left[\left(-k\ln n+\frac{k^2}{2}\ln 2+\frac{17k^2\ln 2}{40}\right)(1+o(1))\right]\\
 &\hspace{5cm}+\sum_{s=s_0+1}^{k}\exp\left[\left(-k\ln n+\frac{k^2}{2}\ln 2+\frac{5k^2\ln 2}{16}\right)(1+o(1))\right]\\
& \leq (1+o(1))\sum_{s=1}^{s_0}\left(\frac{k^2}{2^{k/2}}\right)^s+\exp\left[-\left(1+o(1)\right)\frac{k^2\ln 2}{20}\right]=o(1).
\end{align*}

\section{Proof of Theorem~\ref{th:main} part 2}
\label{proof:2}

The proof strategy is very similar. The difference is that, instead of finding a monotone isomorphism-closed family of graphs such that the respective expected number is asymptotically a constant, we find a set $S_1$ in one of the random graphs (say, $G_1$) that is `non-extendible' in some sense, and show that in the other random graph whp any set $S_2$ that induces a subgraph isomorphic to $G_1[S]$ is extendible.\\

Consider a uniformly distributed random graph $G_n$ on $[n]$, and let $X_k$ be the total number of $k$-sets $S\subset[n]$ such that there is no vertex $x$ outside $S$ adjacent to all vertices in $S$. Then
$$
 {\sf E}X_k={n\choose k}\left(1-\frac{1}{2^k}\right)^{n-k}=e^{f(k)},\quad\text{where}\quad
 f(k)= k\ln  n-k\ln k+k-n2^{-k}+O(\ln k).
$$
Let $k\in\mathbb{R}$ be such that $f(k)=0$. Note that ${\sf E}X_{\lfloor k \rfloor}\leq 1$  and ${\sf E}X_{\lceil k \rceil}\geq 1$. Let $\Delta:=\lceil k\rceil-k$. Then
$$
f(\lceil k\rceil)-
f(k)=k^2\ln 2\left(1-2^{-\Delta}+o(1)\right).
$$
Note that we need ${\sf E}X_{\lceil k \rceil}\to\infty$ in order to find in the first random graph a non-extendible $\lceil k\rceil$-set $S$. Also, we need $\lceil k \rceil!2^{-\lceil k \rceil\choose 2}{\sf E}X_{\lceil k \rceil}\to 0$ in order to show that in the second random graph there is no non-extendible $\lceil k\rceil$-set that induces a graph isomorphic to the graph induced by $S$ in the first random graph.

In order to be certain that ${\sf E}X_{\lceil k \rceil}\to\infty$, we need $1-2^{-\Delta}$ to be bounded away from 0. Moreover, in order to have $\lceil k \rceil!2^{-\lceil k \rceil\choose 2}{\sf E}X_{\lceil k \rceil}\to 0$, we need $1-2^{-\Delta}$ not to be very large. Let us first consider the sufficient restriction $0.05\leq \Delta\leq 0.95$. For other values of $\Delta$, we will slightly modify the sentence in the same way as it was done in Section~\ref{sc:modified_extension}. 

Since $1-2^{-\Delta}$ is bounded away both from $0$ and $1/2$, we get that ${\sf E}X_{\lceil k \rceil}\to\infty$ and $\lceil k \rceil!2^{-\lceil k \rceil\choose 2}{\sf E}X_{\lceil k \rceil}\to 0$. 

On the one hand, it means that ${\sf P}(X_{\lceil k\rceil}>0)\to 1$. Indeed, to show this, it is sufficient (by Chebyshev's inequality) to prove that $\mathrm{\sf E}X_{\lceil k\rceil}(X_{\lceil k\rceil}-1)=(1+o(1))({\sf E}X_{\lceil k\rceil})^2$. Note that this second factorial moment equals to the summation over all choices of pairs of distinct $\lceil k\rceil$-set $S_i,S_j$ of the probability $P_{i,j}$ that each of the sets does not have a vertex outside adjacent to all vertices inside. In Section~\ref{sc:momentsX}, we showed that the contribution of intersecting pairs is $o(({\sf E}X_{\lceil k\rceil})^2)$. The contribution of non-intersection pairs is at most
$$
 {n\choose \lceil k\rceil}^2\left(1-\frac{2}{2^{\lceil k\rceil}}\right)^{n-2\lceil k\rceil}\sim({\sf E}X_{\lceil k\rceil})^2
$$
as needed.

On the other hand, for a given graph $F$ on $[k]$, by the union bound, the probability that there is a set $S$ in $G(n,p)$ that induces a subgraph isomorphic to $F$ such that every vertex outside $S$ is not adjacent to at least 1 vertex of $S$ is at most
$$
 {n\choose \lceil k\rceil}\frac{\lceil k\rceil!}{2^{{\lceil k\rceil\choose 2}}}\left(1-\frac{1}{2^{\lceil k\rceil}}\right)^{n-\lceil k\rceil}\to 0.
$$

Note that $k=\hat k+o(1)$ (recall that $\hat k$ is defined in Section~\ref{sc:intro} as the unique root of the equation $f^{\mathrm{ext}}(k)=0$) implying (since $\Delta\geq 0.05$) that $\lceil k\rceil=\lceil \hat k \rceil=\lfloor\hat k\rfloor +1$. The sentence requires $\lceil k\rceil+1\leq\lfloor\hat k\rfloor +2$ variables.\\

Finally, let either $\Delta<0.05$, or $\Delta>0.95$. Let $Y$ be the total number of tuples $(v_1,v_2,S)$, where $v_1\neq v_2\in[n]$, $S\in{[n]\setminus\{v_1,v_2\}\choose k}$, such that there are no vertices $x$ outside $S\sqcup\{v_1,v_2\}$ adjacent to all vertices from $S$ and at least 1 vertex among $v_1,v_2$. We have
$$
 {\sf E}Y_k={n\choose 2}{n-2\choose k}\left(1-\frac{3}{2^{k+2}}\right)^{n-k-2}=e^{g(k)},
$$
$$ 
 g(k)= (k+2)\ln  n-k\ln k+k-3n2^{-k-2}+O(\ln k).
$$
Let $k'\in\mathbb{R}$ be such that $g(k')=0$. Note that  $k'=k+\log_2 3-2+o(1)$ and
$$
\Delta':=k'-\lfloor k'\rfloor=\Delta+\log_2\frac{3}{2}+o(1)-\left\lfloor\Delta+\log_2\frac{3}{2}+o(1)\right\rfloor\in(0.5,0.65).
$$
Then ${\sf E}Y_{\lceil k'\rceil}\to\infty$ and $\lceil k' \rceil!2^{-\lceil k' \rceil\choose 2}{\sf E}Y_{\lceil k' \rceil}\to 0$ as needed. In the same way as above, applying the bounds from Section~\ref{sc:momentsY} and Chebyshev's inequality, we get that ${\sf P}(Y_{\lceil k'\rceil}>0)\to 1$. Also, for a given graph $F$ on $[\lceil k'\rceil]$, by the union bound, the probability that there are a set $S$ and a pair of distinct vertices $v_1,v_2$ outside $S$ in $G(n,p)$ such that $S$ induces a subgraph isomorphic to $F$ and every common neighbor $u\notin\{v_1,v_2\}$ of $S$ is not adjacent to any of $v_1,v_2$ is at most
$$
 {n\choose 2}{n-2\choose \lceil k\rceil}\frac{\lceil k\rceil!}{2^{{\lceil k\rceil\choose 2}}}\left(1-\frac{3}{2^{\lceil k+2\rceil}}\right)^{n-\lceil k\rceil-2}\to 0.
$$
This finishes the proof of part 2 of Theorem~\ref{th:main} since the number of variables we need is $\lceil k'\rceil+3=\lceil \hat k+\log_2 3-2+o(1)\rceil+3=\lceil \hat k \rceil+2\leq\lfloor \hat k\rfloor +3$.\\

\section{Other logics and models of random graphs}

Let us note that our results are immediately extended to infinitary FO logic with finitely many variables~\cite[Section 11.1]{Libkin} since two graphs that satisfy the extension axiom $\varphi_k$ are not distinguishable in the fragment of this logic with $k+1$ variables as well~\cite[Proposition 11.8]{Libkin}. On the other hand, there are monadic second order sentences such that their truth probability does not converge on $G_n$ to 0 or 1. Actually, 1 {\it existential} monadic quantifier and 2 FO variables is enough~\cite{Popova} to distinguish between two independent uniformly distributed random graph on $[n]$ at least for some sequences of $n$ (note that in the cited paper it is proven that there exist a sentence $\varphi$ and two sequences $n_i,m_i$ such that ${\sf P}(G_{n_i}\models\varphi)\to 0$ and ${\sf P}(G_{m_i}\models\varphi)$ is bounded away from 0; however, it is easy to see that the same method can be applied to show that there exists a sequence $n_i$ such that ${\sf P}(G_{n_i}\models\varphi)$ is bounded away both from 0 and 1 for the same sentence $\varphi$). It would be interesting to know if the same number of variables is sufficient for all large enough $n$ to distinguish between two independent random graphs on $[n]$.

It would be also of interest to study distinguishability of random graphs with other probability distributions. In particular, the classical 0-1 law~\cite{Fagin,Glebskii} can be generalized to the binomial random graph $G(n,p)$ for a wide range of edge probabilities $p$ (note that $p=1/2$ corresponds to the considered in the paper uniform distribution): if $np^{\alpha}\to\infty$ for all $\alpha>0$, then $G(n,p)$ obeys FO 0-1 law~\cite{Spencer}. Our results can be also generalized to the same range of $p=p(n)$. In particular, for every constant $p\in(0,1)$, we will also get concentration of the minimum number of variables required to distinguish between two random graphs in a set of constantly many numbers. Note that, in~\cite{KPSS}, it is proven that, for $p=1/2$, there exists a sequence $\{n_i\}_{i\in\mathbb{N}}$ such that whp the minimum quantifier depth (and, therefore, the number of variables) of a FO sentence defining $G(n_i,1/2)$ is concentrated in a set of 5 consecutive numbers. However, it is not clear if there is a constant-size concentration interval for the whole sequence of $n$ and also for other constant $p\in(0,1)$.

The case $p=n^{-\alpha}$ is different. When $\alpha$ is rational and belongs to $(0,1]$ or $\alpha=1+\frac{1}{m}$ for some $m\in\mathbb{N}$, then $G(n,p)$ does not obey FO 0-1 law~\cite{Shelah} since, for any such $\alpha$, there exists a constant size strictly balanced graph with density $1/\alpha$~\cite{RV}, and the number of isomorphic copies of such a graph in $G(n,n^{-\alpha})$ converges in distribution to a Poisson random variable (see~\cite[Theorem 3.19]{Janson}). Therefore, with probability bounded away from 0, two independent copies of $G(n,p)$ can be distinguished by a FO sentence with constant number of variables (equal to the number of vertices in the mentioned strictly balanced graph). It is also not hard to see that, for every $\varepsilon>0$, there exists a constant $k=k(\varepsilon)$ such that two independent random graphs are $k$-distinguishable with probability at least $1-\varepsilon$, since the joint distribution of the numbers of copies of different strictly balanced graphs with the same density $1/\alpha$ converges to a vector of independent Poisson random variables (see~\cite[Remark 3.20]{Janson}). It is not so clear how to construct a fixed sentence with constantly many variables that distinguishes between two independent random graphs with probability arbitrarily close to $\frac{1}{4}$. However, we believe that this is possible and boils down to constructing strictly balanced graphs with fixed densities and linearly (in the number of vertices) growing number of automorphisms.

If $\alpha>1$ and $\alpha\neq 1+\frac{1}{m}$ for any $m\in\mathbb{N}$, then whp $G(n,p=n^{-\alpha})$ is a forest such that every tree with at most $\lfloor\frac{1}{\alpha}-1\rfloor$ vertices appears as a component in $G(n,p)$ at least $n^{\epsilon}$ times (for sufficiently small constant $\epsilon>0$), and there are not trees of size $\lceil\frac{1}{\alpha}-1\rceil$ (see, e.g.,~\cite[Theorems 3.4, 3.29]{Janson}). Moreover, the number of components isomorphic to a fixed small tree satisfies the central limit theorem \cite{Rucinski}. From this, it is not hard to derive that $n^{\epsilon}$ variables is needed to distinguish between two independent copies of $G(n,p)$, and that this number is not concentrated in any bounded interval whp.

The most intriguing case is when $\alpha\in(0,1)$ is irrational. Then the random graph $G(n,p=n^{-\alpha})$ obeys FO 0-1 law~\cite{Shelah}. So, any fixed sentence does not distinguish between two independent random graphs whp. In the same way as in the case of rational $\alpha$, one may consider an infinite sequence of strictly balanced graphs $H_i$, $i\in\mathbb{N}$, such that the density $\rho(H_i)=\frac{|E(H_i)|}{|V(H_i)|}$ approaches $\frac{1}{\alpha}$ as $i\to\infty$. However, in order to make the expected number of $H_i$ to be asymptotically a constant, we need $|V(H_i)|=\Omega(\ln n)$. Nevertheless, we {\it conjecture} that for every $\varepsilon>0$ and every $\frac{\varepsilon}{n}<p<n^{-\varepsilon}$, the minimum number of variables of a FO sentence that distinguishes between two independent copies of $G(n,p)$ is $o(\ln n)$.

Finally, let us mention other random graph models, in particular random $d$-regular graphs $G_{n,d}$. Note that, in the dense case (say, $d=\frac{1}{2}n(1+o(1))$), FO 0-1 law holds~\cite{Haber}. However, our results can not be immediately translated using, say, sandwiching techniques~\cite{Gao} since we need the coupling $G(n,p_1)\subset G_{n,d}\subset G(n,p_2)$ with $p_1,p_2$ very close to $d$ in order to translate the results, but this is not likely to be true. However, we believe that an analogue of Theorem~\ref{th:main} holds true for $G_{n,d}$ as well. Let us also mention the sparse case of constant degree $d$: two random $d$-regular graphs can be distinguished with constantly many variables since there is no FO 0-1 law due to the convergence to Poisson distribution of the number of constant size cycles~\cite{Wormald}. On the other hand, for sufficiently small $\varepsilon>0$, it is possible to construct 3-regular non-isomorphic graphs on $n$ vertices that can not be distinguished by at most $\varepsilon n$ variables in the first order logic with counting quantifiers (that is even stronger)~\cite{CFI}. It is natural to ask the same question about Cayley graphs. Do there exist 3-regular non-isomorphic Cayley graphs on $n$ vertices that can not be distinguished by at most $\varepsilon n$ variables in the first order logic with counting quantifiers? Note that whp 4 variables is enough to distinguish between two independent uniformly distributed random graphs by a FO sentence with counting~\cite{Babai}. It is also worth noting that distinguishing Cayley graphs in FO logic with counting was studied: for example, in~\cite{Ponomarenko,VZ}, this question is addressed in circulant graphs; see also~\cite{BS}, where the effectiveness of FO logic with counting for distinguishing finite groups is investigated.

\bibliographystyle{amsplain}

\begin{thebibliography}{99}


\bibitem{Babai} L. Babai, L. Ku\v{c}era, {\it Canonical labeling of graphs in linear average time}, Proc. 20th IEEE Symp. Found. Computer Sci., 1979, pp. 39--46.

\bibitem{Billingsley} P. Billingsley, {\it Probability and Measure}, Third Edition, J. Wiley \& Sons, 1995.

\bibitem{BS} J. Brachter, P. Schweitzer, {\it On the Weisfeiler-Leman dimension of finite groups}, Proc. 35th Annual ACM/IEEE Symposium on Logic in Computer Science (LICS), 2020,  pp. 287--300.

\bibitem{CFI} J.-Y. Cai, M. F\"{u}rer, N. Immerman, {\it An optimal lower bound on the number of variables for graph identification}, Combinatorica {\bf 12} (1992) 389--410.

\bibitem{Symmetry} P. Erd\H os and A. R\'{e}nyi, {\it Asymmetric graphs}, Acta Mathematica Hungarica, {\bf 14}:3-4 (1963) 295--315.

\bibitem{Fagin} R. Fagin, \emph{Probabilities in finite models}, J. Symbolic Logic {\bf 41} (1976) 50--58.

\bibitem{Gao} P. Gao, M. Isaev, B. D. McKay, {\it Sandwiching random regular graphs between binomial random graphs}, SIAM Conference on Discrete Algorithms (SODA 2020) 690--701.

\bibitem{Glebskii} Y.V. Glebskii, D.I. Kogan, M.I. Liogon'kii, V.A. Talanov, {\it Range and degree of realizability of formulas the restricted predicate calculus}, Cybernetics {\bf 5} (1969) 142--154 (Russian original: Kibernetica {\bf 2}, 17--27).

\bibitem{Haber} S. Haber, M. Krivelevich, {\it The logic of random regular graphs},  Journal of combinatorics {\bf 1}:4 (2010) 389--440.

\bibitem{IRZZ} M. Isaev, I. Rodionov, R. Zhang, M. Zhukovskii, {\it Extremal independence in discrete random systems}, arXiv:2105.04917.

\bibitem{Janson} S. Janson, T. \L uczak, A. Rucinski, {\it Random graphs}, J. Wiley \& Sons, 2000.

\bibitem{KPSS} J.H. Kim, O. Pikhurko, J.H. Spencer, O. Verbitsky, {\it How complex are random graphs in first order logic?}, Random structures \& algorithms {\bf 26}:1-2 (2005) 119--145. 

\bibitem{Libkin} L. Libkin, {\it Elements of finite model theory}, Springer, 2004. 

\bibitem{Ponomarenko}  I. Ponomarenko, {\it On the WL-dimension of circulant graphs of prime power order}, Algebraic Combinatorics {\bf 6}:6 (2023) 1469--1490.

\bibitem{Popova} S. Popova, M. Zhukovskii, {\it Existential monadic second order logic of undirected graphs: The Le Bars conjecture is false}, Annals of Pure and Applied Logic {\bf 170}:4 (2019) 505--514.

\bibitem{RV} A. Ruci\'{n}ski, A. Vince, {\it Balanced graphs and the problem of subgraphs of a random graph}, Congressus Numerantium, {\bf 49} (1985) 181--190.

\bibitem{Rucinski} A. Ruci\'{n}ski, {\it  When are small subgraphs of a random graph normally distributed?}, Probability theory and related fields, {\bf 78}:1 (1988) 1--10.

\bibitem{Shelah} S. Shelah, J.H. Spencer, {\it Zero-one laws for sparse random graphs}, J. Amer. Math. Soc., \textbf{1} (1988) 97--115.

\bibitem{Spencer} J. H. Spencer, {\it Threshold spectra via the Ehrenfeucht game}, Discrete Applied Math. {\bf 30} (1991)~235--252.

\bibitem{VZ} O. Verbitsky, M. Zhukovskii, {\it Canonization of a random circulant graph by counting walks}, WALCOM 2024. Lecture Notes in Computer Science, vol 14549.

\bibitem{Wormald} N. Wormald, {\it The asymptotic distribution of short cycles in random regular graphs}, J. Combin. Theory Ser. B {\bf 31} (1981) 168--182.

\end{thebibliography}


\begin{dajauthors}
\begin{authorinfo}[itai]
  Itai Benjamini\\
  Weizmann Institute of Science, Department of Mathematics\\
  Rehovot, Israel\\
  itai\imagedot{}benjamini\imageat{}weizmann\imagedot{}ac\imagedot{}il \\
\end{authorinfo}
\begin{authorinfo}[maks]
  Maksim Zhukovskii\\
  The University of Sheffield,  School of Computer Science\\
  Sheffield, UK\\
  m\imagedot{}zhukovskii\imageat{}sheffield\imagedot{}ac\imagedot{}uk \\
\end{authorinfo}
\end{dajauthors}

\end{document}